\documentclass[11pt]{article}

\title{\begin{center}
Singular 
Welschinger invariants
\end{center}}

\author{Eugenii Shustin\thanks{School of Mathematical Sciences, Tel Aviv University,
Ramat Aviv, 69978 Tel Aviv, Israel. E-mail: shustin@tauex.tau.ac.il}}

\date{}

\usepackage{fancyhdr}
\usepackage[all]{xy}
\usepackage{graphicx}
\usepackage{amsmath}
\usepackage{amsthm}
\usepackage{amssymb}

\usepackage{latexsym,amsfonts,amscd,yhmath,epsfig,epic}

\newcommand{\proofend}{\hfill$\Box$\bigskip}

\newcommand{\Q}{\mathbb{Q}}
\newcommand{\R}{\mathbb{R}}
\newcommand{\C}{\mathbb{C}}

\newcommand{\Sing}{\operatorname{Sing}}

\newcommand{\mt}{\operatorname{mult}}
\newcommand{\ord}{\operatorname{ord}}

\DeclareMathOperator{\codim}{codim} %
\newtheorem{theorem}{Theorem}[section]
\newtheorem{proposition}[theorem]{Proposition}
\newtheorem{definition}[theorem]{Definition}

\newtheorem{lemma}[theorem]{Lemma}

\newtheorem{remark}[theorem]{Remark}

\begin{document}

\maketitle

\begin{abstract}
We suggest an invariant way to enumerate nodal and nodal-cuspidal real deformations of real plane curve singularities. The key idea is to assign
Welschinger signs to the counted deformations. Our invariants can be viewed as a local version of Welschinger invariants enumerating
real plane rational curves.
\end{abstract}


\section*{Introduction}

Gromov-Witten invariants of the plane can be identified with the degrees of Severi varieties, which parameterize irreducible plane curves of given degree and genus. As a local version, one can consider a versal deformation
of an isolated plane curve singularity $(C,z)\subset\C^2$ with base $B(C,z)\simeq(\C^n,0)$, 
and the following strata in $B(C,z)$:
\begin{equation}EG^i_{C,z},\quad 1\le i\le\delta(C,z)\ ,\label{le1}\end{equation}
parameterizing deformations with the total $\delta$-invariant greater or equal to $i$;
\begin{equation}EC^k_{C,z},\quad0\le k\le\kappa(C,z)-2\delta(C,z)\ ,\label{le3}\end{equation} parameterizing deformations with the total $\delta$-invariant equal to $\delta(C,z)$ and the total $\kappa$-invariant equal to $2\delta(C,z)+k$ (a necessary information on $\delta$- and $\kappa$-invariants can be found in
\cite{DH} or \cite[Section 3.4]{GLS}). Note also that $EC^0_{C,z}=EG^{\delta(C,z)}_{C,z}$.

The strata (\ref{le1}) are called {\it Severi loci}; among them, ${\mathcal D}_{C,z}:=EG^1_{C,z}$ is the discriminant hypersurface in $B(C,z)$,
and $EG_{C,z}:=EG^{\delta(C,z)}_{C,z}$ is the so-called {\it equigeneric locus}. We call the strata (\ref{le3}) {\it generalized equiclassical loci}, and among them $EC_{C,z}:=EC^{\kappa(C,z)-2\delta(C,z)}_{C,z}$ is the so-called {\it equiclassical locus}. The incidence relations are as follows:
$$EG^i_{C,z}\subsetneq EG^{i+1}_{C,z}\subsetneq EC^k_{C,z}\subsetneq EC^{k+1}_{C,z}$$
for all $1\le i<\delta(C,z)$ and $1\le k<\kappa(C,z)-2\delta(C,z)$.
All these loci are pure-dimensional germs of complex spaces
(cf. \cite{Sch,She}).

A natural problem is to compute the multiplicities of $EG^i_{C,z},EC^k_{C,z}$ for all $i,k$\ \footnote{We understand the multiplicity of a point of an algebraic variety embedded into an affine space as the intersection number at this point
with a generic smooth germ of the complementary dimension (cf. \cite[Chapter 5, Definition 5.9]{Mum}).}. This problem was solved for the equigeneric stratum $EG_{C,z}$ in
\cite{FGS}. In the particular case of an irreducible germ with one Puiseux pair, i.e., topologically equivalent to $x^p+y^q=0$, $2\le p<q$, $\gcd(p,q)=1$, one has (see \cite[Proposition 4.3]{Be} and \cite[Section G]{FGS})
$$\mt EG_{C,z}=\frac{1}{p+q}\binom{p+q}{p}\ .$$
The multiplicities of
all Severi loci $EG^i_{C,z}$ were expressed in \cite{She} in terms of the Euler characteristics of Hilbert schemes of points on curve germs representing a given singularity.
The multiplicities of the equiclassical loci $EC^k_{C,z}$ are not known except for the case of the smoothness mentioned in \cite[Theorems 2 and 27]{Di}.

The multiplicity admits an enumerative interpretation: it can be regarded as the number of
intersection points of a locus $V\subset B(C,z)$ with a generic affine subspace $L\subset B(C,z)$ of the complementary dimension (equal to $\codim V$) chosen to be transversal to the tangent cone $\widehat T_0V$.

The {\bf goal} of this note is to define {\it real multiplicities} of the Severi loci (\ref{le1}) and of the generalized equiclassical loci (\ref{le3}). Let the singularity $(C,z)$ be real\footnote{Under the {\it real} object we always understand a complex object invariant with respect to the complex conjugation.} Then the Severi loci and the generalized equiclassical loci are defined over the reals. Thus, given such a locus $V$, we count real intersection points of $V$ with a
generic real affine subspace $L\subset B(C,z)$ of the complementary dimension. Our {\bf main result} is that, in certain cases, the count of real intersection points of $V$ and $L$ equipped with Welschinger-type signs is invariant, i.e., does not depend on the choice of $L$.
We were motivated by
\cite[Lemma 15]{IKS2}, which, in fact, states the existence of a Welschinger type invariant for
the equigeneric stratum $EG_{C,z}$. In this note, we go further and prove the existence of similar Welschinger type invariants for $EG^{\delta(C,z)-1}_{C,z}$ (see Proposition \ref{lp2} in Section \ref{lsec1}) and for $EG^1_{C,z}={\mathcal D}_{C,z}\subset B(C,z)$ (see Proposition \ref{lp4} in Section \ref{lsec1})
as well as for all the loci $EC^k_{C,z}$ (see Proposition \ref{lp3} in Section \ref{lsec2}).

We remark that a similar
enumeration of real plane rational curves with at least one cusp is not invariant, i.e., depends on the choice of point constraints (cf. \cite{Wel}).

As an example, we perform computations for singularities of type $A_n$ (see Section \ref{lsec4}).

\smallskip

{\bf Acknowledgements.} The work at this paper has been supported by the Israeli Science Foundation grants no. 176/15 and 501/18, and by the Bauer-Neuman chair in Real and Complex Geometry. I also would like to thank Stephan Snegirov, with whom I discussed the computational part of the work.

\section{Singular Welschinger numbers}\label{sec-pr}

We shortly recall definitions and basic properties of objects of our interest. Details can be found in \cite{DH} and
\cite[Chaper II]{GLS}.

Let $(C,z)$ be the germ of a plane complex analytic curve $C$ at its isolated singular point $z=(0,0)\in\C^2$, which is given by an analytic equation $f(x,y)=0$, $f\in\C\{x,y\}$.
We shortly call it {\it singularity}. The Milnor ball $D(C,z)\subset\C^2$ is a closed ball centered at $z$
such that $C\cap D(C,z)$ is closed and smooth outside $z$ with the boundary $\partial (C\cap D(C,z))\subset\partial D(C,z)$, and
the intersection of $C$ with any $3$-sphere in $D(C,z)$ centered at $z$ is transversal.
Pick integer $N>0$ and consider the small neighborhood $B(C,z)$ of $0$ in the space (which is a $\C$-algebra)
$R(C,z):=\C\{x,y\}/(\langle f\rangle+{\mathfrak m}_z^N)$, where ${\mathfrak m}_z\subset\C\{x,y\}$ is the maximal ideal.
We can suppose that, for any $\varphi\in B(C,z)$, the curve ${\mathcal C}_\varphi:=\{f+\varphi=0\}\cap D(C,z)$ has only isolated singularities in $D(C,z)$, is smooth along $\partial D(C,z)$, and intersects the sphere $\partial D(C,z)$ transversally.
It is well-known that
the deformation $\pi:{\mathcal C}\to B(C,z)$ of $(C,z)$, where $\pi^{-1}(\varphi)={\mathcal C}_\varphi$, is versal for $N>0$ sufficiently large (cf. \cite[Page 165]{Ar} or \cite[Section 3]{DH}).
The space $B(C,z)$ contains the equigeneric stratum $EG_{C,z}\subset B(C,z)$, formed by $\varphi\in B(C,z)$ such that ${\mathcal C}_\varphi$ has the total $\delta$-invariant equal to
$\delta(C,z)$ (the maximal possible value), the equiclassical locus $EC_{C,z}\subset
EG_{C,z}\subset B(C,z)$, formed by $\varphi\in EG_{C,z}$ such that ${\mathcal C}_\varphi$ has the total $\kappa$-invariant
equal to $\kappa(C,z)$ (also the maximal possible value), and the discriminant
$${\mathcal D}_{C,z}=\{\varphi\in B(C,z)\ :\ {\mathcal C}_\varphi\ \text{is singular}\}\ .$$

The following statement summarizes some known facts on the above strata (see \cite[Theorems 1.1, 1.3, 4.15, 4.17, 5.5, Corollary 5.13]{DH} and \cite[Theorems 2 and 27]{Di}).

\begin{lemma}\label{ll1}
(1) The stratum $EG_{C,z}$ is irreducible of codimension $\delta(C,z)$; it is smooth iff all irreducible components of
$(C,z)$ (which we call local branches of $(C,z)$) are smooth; in general, the normalization of $EG_{C,z}$ is smooth and projects one-to-one onto $EG_{C,z}$. The tangent cone $\widehat T_0EG_{C,z}$ is the linear space $J^{cond}_{C,z}/{\mathfrak m}_z^N$ of codimension
$\delta(C,z)$, where $J^{cond}_{C,z}\subset\C\{x,y\}/\langle f\rangle$ is the conductor ideal. Furthermore, $EG_{C,z}$ contains an open dense subset
$EG^*_{C,z}$ that parameterizes the curves ${\mathcal C}_\varphi$ having $\delta(C,z)$ nodes as their only singularities.

(2) The stratum $EC_{C,z}$ is irreducible of codimension $\kappa(C,z)-\delta(C,z)$; it is smooth iff each local branch of
$(C,z)$ either is smooth, or has topological type $x^m+y^{m+1}=0$ with $m\ge 2$;
in general, the normalization of $EC_{C,z}$ is smooth and projects one-to-one onto $EC_{C,z}$. The tangent cone $\widehat T_0EC_{C,z}$ is the linear space $J^{ec}_{C,z}/{\mathfrak m}_z^N$ of codimension
$\kappa(C,z)-\delta(C,z)$, where $J^{ec}_{C,z}\subset\C\{x,y\}/\langle f\rangle$ is the equiclassical ideal.
Furthermore, the stratum $EC_{C,z}$ contains an open dense subset $EC^*_{C,z}$ that parameterizes the curves
${\mathcal C}_\varphi$ having $3\delta(C,z)-\kappa(C,z)$ nodes and $\kappa(C,z)-2\delta(C,z)$ ordinary cusps as their only singularities.

(3) The discriminant ${\mathcal D}_{C,z}$ is an irreducible hypersurface with the tangent cone
$\widehat T_0{\mathcal D}_{C,z}={\mathfrak m}_z/(\langle f\rangle+{\mathfrak m}_z^N)$. Furthermore, an open dense subset
$\quad$ ${\mathcal D}^*_{C,z}\subset{\mathcal D}_{C,z}$ parameterizes the curves ${\mathcal C}_\varphi$ having one node and no other singularities.
\end{lemma}

In the same way one can establish similar properties of the Severi loci (\ref{le1}) and generalized equiclassical loci (\ref{le3}).

\begin{lemma}\label{ll4}
(1) Each Severi locus $EG^i_{C,z}$ is a (possibly reducible) germ of a complex space of pure codimension $i$. A generic element of each component of $EG^i_{C,z}$ is a curve with $i$ nodes as its only singularities.

(2) Each generalized equiclassical locus $EC^k_{C,z}$ is a (possibly reducible) germ of a complex space of pure codimension $\delta(C,z)+k$. A generic element of each component of $EC^k_{C,z}$ is a curve with $\delta(C,z)-k$ nodes and $k$ ordinary cusps as its only singularities.
\end{lemma}

It is well-known that $\mt{\mathcal D}_{C,z}=\mu(C,z)$ (the Milnor number), $\mt EG_{C,z}$ has been computed in
\cite{FGS} as the Euler characteristic of an appropriate compactified Jacobian.

Now we switch to the real setting. We call the complex space $V$ real if it is invariant under the (natural) action of the complex conjugation and denote by $\R V$ its real point set. Suppose that $(C,z)$ is real.

\begin{definition}\label{ld1}
Let $V\subset B(C,z)$ be an equivariant union of irreducible components of either
a Severi locus $EG^i_{C,z}$, $1\le i\le\delta(C,z)$, or a generalized equiclassical locus $EC^k_{C,z}$,
$1\le k\le\kappa(C,z)-2\delta(C,z)$, and let $\widehat T_0V$ be a linear subspace of $R(C,z)$ of dimension $\dim V$. Assume that $L_0\subset R(C,z)$ is a real linear subspace of dimension
$\dim L_0=\codim_{B(C,z)}V$, which meets $\widehat T_0V$ only at the origin, and let $U(L_0)$ be a neighborhood
of the origin such that $L_0\cap V\cap U(L_0)=\{0\}$. For a sufficiently close to $L_0$ real affine space $L\subset R(C,z)$ of dimension $\dim L=\dim L_0$, intersecting $V\cap U(L_0)$ along $V^*$ and
with total multiplicity
$\mt V$, we set
$$W(C,z,V,L)=\sum_{\varphi\in L\cap\R V\cap U(L_0)}w(\varphi),\quad\text{where}\ w(\varphi)=(-1)^{s(\varphi)
+ic(\varphi)}\ ,$$ with
$s(\varphi)$ being the number of real elliptic\footnote{A real node is called elliptic if it is equivariantly isomorphic to $x^2+y^2=0$.} nodes of ${\mathcal C}_\varphi$ and $ic(\varphi$ the number of pairs of complex conjugate cusps of ${\mathcal C}_\varphi$. In case of $V=EG_{C,z}$ or
$EC_{C,z}$, we write $W^{eg}(C,z,L)$ or $W^{ec}(C,z,L)$, respectively. \end{definition}

In what follows we examine the dependence on $L$ and prove some invariance statements.

\section{Singular Welschinger invariant $W^{eg}(C,z)$}
The following statement is a consequence of \cite[Lemma 15]{IKS2}. We provide a proof, since in a similar manner we treat other instances of the invariance.

\begin{proposition}\label{lp1}
Given a real singularity $(C,z)$, the number $W^{eg}(C,z,L)$ does not depend on the choice of $L$.
\end{proposition}

{\bf Proof.}
Let $L'_0,L''_0\subset R(C,z)$ be two real linear subspaces of dimension $\delta(C,z)$ transversally intersecting $T_0EG_{C,z}$ at the origin, and let $L',L''\subset R(C,z)$ be real affine subspaces of dimension
$\delta(C,z)$, which are sufficiently close to $L'_0,L''_0$, respectively, in the sense of Definition \ref{ld1}. We can connect the pairs $(L'_0,L')$ and $(L''_0,L'')$ by a generic smooth path
$\{L_0(t),L(t)\}_{t\in0,1]}$ consisting of real linear subspaces $L_0(t)$ of
$R(C,z)$ of dimension $\delta(C,z)$, which are transversal to $T_0EG_{C,z}$, and real affine subspaces $L(t)$ of dimension
$\delta(C,z)$ sufficiently close to $L_0(t)$ in the sense of Definition \ref{ld1}, $0\le t\le 1$.
It follows from Lemma \ref{ll1}(1) that, for all $t\in[0,1]$, the space $L(t)$ intersects $EG_{C,z}$
transversally at each element of $L(t)\cap EG_{C,z}$.
Furthermore, all but finitely many spaces $L(t)$ intersect $EG_{C,z}$ along $EG^*_{C,z}$, transversally at each intersection point. The remaining finite subset $F\subset(0,1)$ is such that, for any $\hat t\in F$, the
intersection $L(\hat t)\cap EG_{C,z}$ consists of elements of $EG^*_{C,z}$ and one real element $\varphi$ belonging to a codimension one substratum of $EG_{C,z}$. The classification of these codimension one substrata is known (see, for instance
\cite[Theorem 1.4]{DH}): an element $\varphi$ of such a substratum is as follows:
\begin{enumerate}\item[(n1)] either ${\mathcal C}_\varphi$ has an ordinary cusp $A_2$ and $\delta(C,z)-1$ nodes,
\item[(n2)] or ${\mathcal C}_\varphi$ has a tacnode $A_3$ and $\delta(C,z)-2$ nodes,
\item[(n3)] or ${\mathcal C}_\varphi$ has a triple point $D_4$ and $\delta(C,z)-3$ nodes.
\end{enumerate}

In cases (n2) and (n3), the stratum $EG_{C,z}$ is smooth at $\varphi$ (cf. Lemma
\ref{ll1}(1)), and the deformation of ${\mathcal C}_\varphi$ under the variation of $L(t)$ induces
independent equivariant deformations of all (smooth) local branches
of ${\mathcal C}_varphi$ at the non-nodal singular point. Then the exponent $s(\psi)$ (see Definition \ref{ld1}) for any real nodal curve ${\mathcal C}_\psi$, $\psi\in EG_{C,z}$ close to ${\mathcal C}_\varphi$ always equals modulo $2$ the number of elliptic nodes of ${\mathcal C}_\varphi$ plus the intersection number of complex conjugate local branches of ${\mathcal C}_\varphi$ at the non-nodal singular point. Thus, the crossing of these strata does not affect $W^{eg}(C,z,L(t))$.

In case (n1), the germ of $B(C,z)$ at $\varphi$ can be represented as $\qquad\qquad$ $B(A_2)\times B(A_1)^{\delta(C,z)-1}\times(\C^{n-\delta(C,z)-1},0)$ (cf. \cite[Proposition I.1.14 and Theorem I.1.15]{GLS}
and \cite[Lemma 13]{IKS2}), where $n=\dim B(C,z)$, $B(A_2)\simeq(\C^2,0)$ is a miniversal deformation base of an ordinary cusp, which we without loss of generality can identify with the base of the deformation
$\{y^2-x^3-\alpha x-\beta\ :\ \alpha,\beta\in(\C^2,0)\}$, and $B(A_1)\simeq(\C,0)$ stands for the versal deformation of an ordinary node. Here $$(EG_{C,z},b)=EG(A_2)\times EG(A_1)^{\delta(C,z)-1}\times(\C^{n-\delta(C,z)-1},0)\ ,$$ where $n=\dim B(C,z)$ and
$$EG(A_2)=\left\{\frac{\alpha^3}{27}-\frac{\beta^2}{4}\right\}\subset B(A_2),\quad EG(A_1)=\{0\}\subset B(A_1)\ ,$$
$$T_bEG_{C,z}=EG(A_2)\times\{0\}^{\delta(C,z)-1}\times\C^{n-\delta(C,z)-1}\ .$$ Then the transversality of the intersection of
$L(\hat t)$ and $T_\varphi EG_{C,z}$ yields that the family $\{L(t)\}_{|t-\hat t|<\eta}$ projects to the family of
smooth curves $\{L^1(t)\}_{|t-\hat t|<\eta}$ transversal to $T_0EG(A_2)=\{\beta=0\}$. It is easy to see that either
$L^1(t)$ does not intersect $EG(A_2)$ in real points, or it intersects $EG(A_2)$ in two real points
$(\alpha_1,\beta_1)$, $(\alpha_2,\beta_2)$ with $\beta_1<0<\beta_2$, where the former point corresponds to a real curve with a hyperbolic node in a neighborhood of the cusp, while the latter one - to a real curve
with an elliptic node. Hence, the Welschinger signs of these intersections of $L^1(t)$ with $EG(A_2)$ cancel out, which confirms the constancy of $W^{eg}(C,z,L(t))$, $|t-\hat t|<\eta$, in the considered wall-crossing.
\proofend

We mention also two more useful properties of the invariant $W^{eg}(C,z)$.

\begin{lemma}\label{ll3}
(1) The number $W^{eg}(C,z)$ is an invariant of a real equisingular deformation class. That is, if $(C_t,z)_{t\in[0,1]}$ is an equisingular\footnote{``Equisingular" means ``preserving the (complex) topological type".} family of real singularities, then $W^{eg}(C_0,z)=W^{eg}(C_1,z)$.

(2) Let $(C,z)=\bigcup_{i}(C_i,z)$ be the decomposition of a real singularity $(C,z)$ into irreducible over $\R$ components. Then $W^{eg}(C,z)=\prod_iW^{eg}(C_i,z)$.
\end{lemma}

{\bf Proof.}
(1) It is sufficient to verify the local constancy of $W^{eg}(C,z)$ in real equisingular deformations. Recall that the equisingular stratum $ES_{C,z}\subset B(C,z)$ is a smooth subvariety germ. Furthermore, for $N>\mu(C,z)+1$, the germ of $B(C,z)$ at any point $\psi\in ES_{C,z}$ is a versal deformation of of the singularity ${\mathcal C}_\psi$. Then the equality $W^{eg}_{C,z}=W^{eg}_{{\mathcal C}_\psi}$ follows by the argument in the proof of Proposition \ref{lp1}.

(2) The second statement of the lemma follows from the fact that a equigeneric deformation of $(C,z)$ induces independent equigeneric deformations of the components $(C_i,z)$ and vice versa (see \cite[Theorem 1, page 73]{Tei}, \cite[Corollary 3.3.1]{ChL}, and also
    \cite[Theorem II.2.56]{GLS}), and from the fact that the deformed components $(C_i,z)$ and $C_j,z)$, $i\ne j$, can intersect only in hyperbolic real nodes and in complex conjugate nodes.
\proofend

\section{Singular Welschinger invariants associated with $EG^{\delta(C,z)-1}_{C,z}$ and
${\mathcal D}_{C,z}$}\label{lsec1}

The key ingredient of the proof of Proposition \ref{lp1} is that the tangent cone to the equigeneric stratum $EG_{C,z}$ is a linear space of dimension equal to $\dim EG_{C,z}$.
We intend to establish a similar statement for $EG^{\delta(C,z)-1}_{C,z}$.

Recall the following fact used in the sequel:
By \cite[Theorem 1.1]{Sch} the closure of each irreducible component of $EG^{\delta(C,z)-1}_{C,z}$ contains $EG_{C,z}$, and a generic element of such a component can be obtained by smoothing a node of an element of $EG^*_{C,z}$.

\begin{lemma}\label{ll2}
For the following real substrata $V\subset EG^{\delta(C,z)-1}_{C,z}$, the tangent cones $\widehat T_0\R V$ are linear subspaces of $\R R(C,z)$ of (real)
codimension $\qquad\qquad\qquad$ \mbox{$\delta(C,z)-1=\codim \R EG^{\delta(C,z)-1}_{C,z}$}:
\begin{enumerate}\item[(i)] $(C,z)$ contains a real singular local branch $(C',z)$, and
$V\subset EG^{\delta(C,z)-1}_{C,z}$ is the union of those irreducible components of $EG^{\delta(C,z)-1}_{C,z}$, which contain nodal curves obtained from the curves ${\mathcal C}_\varphi$, $\varphi\in EG^*_{C,z}$, by smoothing out a real node on the component of ${\mathcal C}_\varphi$ corresponding to the local branch $(C',z)$;
\item[(ii)] $(C,z)$ contains a pair of complex conjugate local branches $(C',z)$, $(C'',z)$, and
$V\subset EG^{\delta(C,z)-1}_{C,z}$ is the union of those irreducible components of $EG^{\delta(C,z)-1}_{C,z}$, which contain nodal curves obtained from the curves ${\mathcal C}_\varphi$, $\varphi\in EG^*_{C,z}$, by smoothing out a real
intersection point on the components of ${\mathcal C}_\varphi$ corresponding to the local branches $(C',z)$, $(C'',z)$.
\end{enumerate}
\end{lemma}

{\bf Proof.}
(i) Notice, first, that $V$ can be identified with $EG^{\delta(C',z)-1}_{C',z}\times EG_{C'',z}$, where
$(C'',z)$ is the union of the local branches of $(C,z)$ different from $(C',z)$. Hence, we can simply assume that $(C,z)$ is irreducible.

If ${\mathcal C}_\varphi$, $\varphi\in EG^{\delta(C,z)-1}_{C,z}$, has precisely $\delta(C,z)-1$ nodes as its only singularities, then
the tangent space $T_\varphi EG^{\delta(C,z)-1}_{C,z}$ can be identified with the space of elements
$\psi\in R(C,z)$ vanishing at the nodes of ${\mathcal C}_{\varphi}$. It has codimension $\delta(C,z)-1$, and we have the following bound for the intersection:
$$({\mathcal C}_{\psi}\cdot{\mathcal C}_\varphi)\ge2\delta(C,z)-2,\quad \psi\in T_\varphi EG^{-1}_{C,z}\ .$$
Hence, any limit of the tangent spaces $T_\varphi EG^{\delta(C,z)-1}_{C,z}$ as $\varphi\to0$ is contained in the linear space
$$\{\psi\in R(C,z)\ :\ \ord\\psi\big|_{(C,z)}\ge2\delta(C,z)-2\}$$ of codimension at most $\delta(C,z)-1$.
By \cite[Proposition 5.8.6]{CA} we have
$$\{\psi\in R(C,z)\ :\ \ord\psi\big|_{(C,z)}\ge2\delta(C,z)-1\}$$
$$\qquad=
\{\psi\in R(C,z)\ :\ \ord\psi\big|_{(C,z)}\ge2\delta(C,z)\}=
J^{cond}_{C,z}/{\mathfrak m}_z^N\ .$$ Hence
$$\codim\{\psi\in R(C,z)\ :\ \ord\psi\big|_{(C,z)}\ge2\delta(C,z)-2\}$$
$$\ge\codim J^{cond}_{C,z}/{\mathfrak m}_z^N-1=\delta(C,z)-1\ ,$$ and we are done.

\smallskip
(ii) As in the preceding case, we can assume that $(C,z)=(C',z)\cup(C'',z)$. The above argument yields that
the limits of the tangent spaces $T_\varphi\R V$ as $\varphi\in\R V^*$ tends to $0$, are contained in the linear subspace
$$\{\psi\in\R R(C,z)\ :\ \ord\psi\big|_{(C',z)}=\ord\psi\big|_{(C'',z)}\ge2\delta(C',z)+(C'\cdot C'')_z-1\}$$
which then must be of (real) codimension at most $\delta(C,z)-1$. So, it remains to show that the latter codimension equals exactly $\delta(C,z)-1$, and we will prove that the complex codimension of the space
$$\{\psi\in R(C,z)\ :\ \ord\psi\big|_{(C',z)}=\ord\psi\big|_{(C'',z)}\ge2\delta(C',z)+(C'\cdot C'')_z-1\}$$
is at least $\delta(C,z)-1$. Namely, we just impose an extra linear condition and show that the resulting space
$$\Lambda=\{\psi\in R(C,z)\ :\ \ord\psi\big|_{(C',z)}\ge2\delta(C',z)+(C'\cdot C'')_z$$
$$\qquad\qquad\qquad\qquad\ord\psi\big|_{(C'',z)}\ge2\delta(C'',z)+(C'\cdot C'')_z-1\}$$ has codimension $\ge\delta(C,z)$.
Write $f=f'f''$, where $f'=0$ and $f''=0$ are equations of $(C',z)$, $(C'',z)$, respectively.
By the Noether's theorem in the form of \cite[Theorem II.2.1.26]{GLS2}, any $\psi\in\Lambda$ can be represented as
$\psi=af'+bf''$, where $a,b\in R(C,z)$ and
$$\ord a\big|_{(C'',z)}\ge2\delta(C'',z)-1,\quad \ord b\big|_{(C',z)}\ge2\delta(C',z)\ .$$ Again by
\cite[Proposition 5.8.6]{CA}, the former inequality yields
$$\ord a\big|_{(C'',z)}\ge2\delta(C'',z)\ ,$$ which finally implies that
$\Lambda\subset J^{cond}_{C,z}/{\mathfrak m}_z^N$, and hence $\codim\Lambda\ge\delta(C,z)$.
\proofend

\begin{proposition}\label{lp2}
Let $V\subset EG^{\delta(C,z)-1}_{C,z}$ satisfy the hypotheses of one of the cases in Lemma \ref{ll2}. Then
$W(C,z,V,L)$ does not depend on the choice of the real affine space $L$ as in Definition \ref{ld1}.
\end{proposition}

{\bf Proof.}
We closely follow the argument in the proof of Proposition \ref{lp1}. The classification of codimension one substrata of $V$ contains the cases (n1)-(n3) as in the proof of proposition \ref{lp1}, and one additional case:
\begin{enumerate}\item[(n4)] the substratum is $EG_{C,z}$ (i.e., its generic element $\varphi$ has $\delta(C,z)$ nodes).
\end{enumerate}
The analysis of the cases (a)-(c) literally coincides with that in the proof of Proposition \ref{lp1}.
In case (d), the germ of $\R V$ at $\varphi$ consists of $k$ pairwise transversal smooth real germs of codimension $\delta(C,z)-1$ in $\R R(C,z)$, where $k$ is the number of such real nodes $p$ of the curve
${\mathcal C}_\varphi$ that the smoothing of $p$ yields an element of $\R V$ (depending on $V$ as defined in Lemma \ref{ll2}). For any smooth germ $M$ in this union, the intersection of $L(t)\cap M$, $0<|t-\hat t|<\eta$, yields a curve ${\mathcal C}_\psi$ whose Welschinger sign depends only on the real nodes of
${\mathcal C}_\varphi$ different from $p$, and hence does not depend on $t$.
\proofend

By Lemma \ref{ll1}(3), the tangent cone $\widehat T_0{\mathcal D}_{C,z}$ is a hyperplane. As in the preceding case, this yields

\begin{proposition}\label{lp4}
Given a real singularity $(C,z)$, the number $\qquad\qquad$ $W^{discr}(C,z,L):=W(C,z,{\mathcal D}_{C,z},L)$ does not depend on the choice of a real line $L$.
\end{proposition}

The proof literally follows the argument in the proof of Propositions \ref{lp1} and \ref{lp2}.

\section{Singular Welschinger invariants associated with $EC^k_{C,z}$}\label{lsec2}

We start with the equiclassical stratum $EC_{C,z}$, which is the most interesting.

\begin{proposition}\label{lp3}
(1) Given a real singularity $(C,z)$, the number $W^{ec}(C,z,L)$ does not depend on the choice of $L$.

(2) The number $W^{ec}(C,z)$ is an invariant of a real equisingular deformation class. That is, if $(C_t,z)_{t\in[0,1]}$ is an equisingular family of real singularities, then $W^{ec}(C_0,z)=W^{ec}(C_1,z)$.

(3) Let $(C,z)=\bigcup_{i}(C_i,z)$ be the decomposition of a real singularity $(C,z)$ into irreducible over $\R$ components. Then $W^{ec}(C,z)=\prod_iW^{ec}(C_i,z)$.
\end{proposition}

{\bf Proof.}
Again the proof follows the argument in the proof of Proposition \ref{lp1}. So, we accept the initial setting and the notations in the proof of Proposition
\ref{lp1}. Then we study the wall-crossings that correspond to codimension one substrata in $EC_{C,z}$. If $\varphi\in EC_{C,z}$ is a general element of a codimension one substratum, then
\begin{enumerate}\item[(n1')] either ${\mathcal C}_\varphi$ has $3\delta(C,z)-\kappa(C,z)-1$ nodes and $\kappa(C,z)-2\delta(C,z)+1$ cusps,
\item[(n2')] or ${\mathcal C}_\varphi$ has $3\delta(C,z)-\kappa(C,z)-2$ nodes, $\kappa(C,z)-2\delta(C,z)$ cusps, and one tacnode $A_3$,
\item[(n3')] or ${\mathcal C}_\varphi$ has $3\delta(C,z)-\kappa(C,z)-3$ nodes, $\kappa(C,z)-2\delta(C,z)$ cusps, and one triple point $D_4$,
\item[(c1')] or ${\mathcal C}_\varphi$ has $3\delta(C,z)-\kappa(C,z)-1$ nodes, $\kappa(C,z)-2\delta(C,z)-1$ cusps, and one singularity $A_4$,
\item[(c2')] or ${\mathcal C}_\varphi$ has $3\delta(C,z)-\kappa(C,z)-2$ nodes, $\kappa(C,z)-2\delta(C,z)-1$ cusps, and one singularity $D_5$,
\item[(c3')] or ${\mathcal C}_\varphi$ has $3\delta(C,z)-\kappa(C,z)-1$ nodes, $\kappa(C,z)-2\delta(C,z)-2$ cusps, and one singularity $E_6$.
\end{enumerate}
First, we notice that the wall-crossings of types (n1'), (n2'), (n3') are completely similar to the wall-crossing (n1), (n2), (n3), respectively, considered in the proof of Proposition \ref{lp1}, since they involve only the nodal part of the singularities of degenerating elements of $\R EC^*_{C,z}$. Hence, the constancy of $W^{ec}(C,z,L(t))$, $|t-t^*|<\eta$, follows in the same way.

Next we explain why (c1'), (c2'), (c3') are the only codimension one substrata of $EC_{C,z}$ that involve cusps of the degenerating elements of $\R EC^*_{C,z}$. To this end, we show that, any other collection of singularities of ${\mathcal C}_\varphi$ can be deformed into $3\delta(C,z)-\kappa(C,z)$ nodes and $\kappa(C,z)-2\delta(C,z)$ cusps in two successive non-equisingular deformations. By our assumption, at least one of the non-nodal-cuspidal singularities of ${\mathcal C}_\varphi$ must contain a singular local branch. Thus,
\begin{itemize}\item if ${\mathcal C}_\varphi$ has at least two non-nodal-cuspidal singularity, we, first, deform one such singularity into nodes and cusps (along its equiclassical deformation), then all other singularities;
\item if the non-nodal-cuspidal singularity of ${\mathcal C}_\varphi$ has at least three local branches (one of which denoted $P$ is singular), we, first, shift away a branch, different from $P$, then equiclassically deform the obtained curve into a nodal-cuspidal one;
    \item if the non-nodal-cuspidal singularity of ${\mathcal C}_\varphi$ has two singular branches $P_1,P_2$, we, first, shift $P_2$ so that $P_2$ remains centered at a smooth point of $P_1$, then equiclassically deform the obtained curve into a nodal-cuspidal one;
        \item if the non-nodal-cuspidal singularity of ${\mathcal C}_\varphi$ has two branches, $P_1$  smooth and $P_2$ singular, which is different from an ordinary cusp, then we, first, equiclassically deform the local branch $P_2$ into nodes and (necessarily appearing) cusps, while keeping one cusp centered on $P_1$, then deform the obtained triple singularity into
            nodes and one cusp;
            \item if the non-nodal-cuspidal singularity of ${\mathcal C}_\varphi$ has two branches, $P_1$  smooth and $P_2$ singular of type $A_2$, which is tangent to $P_1$, then we, first, rotate $P_1$ so that it becomes transversal to $P_2$, then deform the obtained singularity $D_5$ into two nodes and one cusp;
                \item if the non-nodal-cuspidal singularity of ${\mathcal C}_\varphi$ is unibranch either of multiplicity $m\ge3$ and not of the topological type $y^m+x^{m+1}=0$, or of multiplicity $2$ and not of type $A_4$, then we, first, equigenerically deform this singularity into some nodes and a singularity of topological type $y^m+x^{m+1}=0$, if $m\ge3$, or a singularity $A_4$, if $m=2$ (this can be done by the blow-up construction
                    as in the proof of \cite[Theorem 1]{AC}, see also \cite[Section 2.1]{LS}), then equiclassically deform the obtained curve into a nodal-cuspidal one;
                    \item if the non-nodal-cuspidal singularity of ${\mathcal C}_\varphi$ is of the topological type $y^m+x^{m+1}=0$, $m\ge4$, then the codimension of the its equisingular stratum in a versal deformation base equals $\frac{m^2+3m}{2}-3$, while the codimension of the equiclassical stratum equal
                        $$\kappa(\{y^m+x^{m+1}=0\})-\delta(\{y^m+x^{m+1}=0\})=\frac{m^2+m}{2}-1$$
                        $$=\left(\frac{m^2+3m}{2}-3\right)-(m-2)\le\left(\frac{m^2+3m}{2}-3\right)-2\ .$$
                    \end{itemize}

Now we analyze the wall-crossings of type (c1'), (c2'), and (c3') as described above.

In case (c1'), the miniversal unfolding of an $A_4$ singularity $y^2=x^5$ is given by the family
$y^2=x^5+a_3x^3+a_2x^2+a_1x+a_0$ with the base $B=\{(a_0,...,a_3)\in(\C^4,0)\}$, while the equiclassical locus
$EC\subset B$ is a curve given by $y^2=(x-2t)^3(x+3t)^2$, $t\in(\C,0)$. This curve has an ordinary cusp at the origin. The natural projection of the germ of $B(C,z)$ at ${\mathcal C}_\varphi$ onto $B$ takes the affine spaces $L(t)$, $|t-t^*|<\eta$, to real three-dimensional affine spaces transversal to the tangent line to
$EC$ at the origin. Similarly to the case (n1) in the proof of Proposition \ref{lp1}, in the considered bifurcation, two real intersections with $EC$, one corresponding to a curve with a cusp and a hyperbolic node and the other corresponding to a curve with a cusp and an elliptic node, turns in the wall-crossing into two complex conjugate intersections, and hence the constancy of $W^{ec}(C,z,L(t))$, $|t-t^*|<\eta$, follows.

In case (c2'), the equiclassical locus in a miniversal deformation base of a singularity $D_5$ given, say, by $x(y^2-x^3)=0$ is smooth
and can be described by a family $(x-t)(y^2-x^3)=0$. So, in the considered wall-crossing a real curve with a cusp and two hyperbolic nodes turns into a curve with a cusp and two complex conjugate nodes,
and hence the constancy of $W^{ec}(C,z,L(t))$, $|t-t^*|<\eta$, follows.

In case (c3'), again the equiclassical locus in a miniversal deformation base of of a singularity $E_6$ 
is smooth (cf. \cite[Theorem 27]{Di}) and one-dimensional. It is not difficult to show that one half branch of $\R EC(E_6)$ parameterizes curves with two real cusps and one hyperbolic node, while the other half branch parameterizes curves with two complex conjugate cusps and one elliptic node. Thus,
the constancy of $W^{ec}(C,z,L(t))$, $|t-t^*|<\eta$, follows.
\proofend

The other loci $EC^k_{C,z}$, $1\le k<\kappa(C,z)-2\delta(C,z)$, may be reducible. Assume that $(C,z)=(C_1,z)\cup...\cup(C_s,z)$ is the splitting into irreducible (over $\C$) components. Given a partition
$\overline k=(k_1,...,k_s)$ such that
\begin{equation}k_1+...+k_s=k,\quad0\le k_i\le\kappa(C_i,z)-2\delta(C_i,z),\ i=1,...,s,
\label{le4}\end{equation}
we define the substratum $EC^{\overline k}_{C,z}\subset EC^k_{C,z}$, which is the union of those irreducible components of $EC^k_{C,z}$ whose generic elements $\varphi$ are such that
${\mathcal C}_\varphi={\mathcal C}_{1,\varphi}\cup...\cup{\mathcal C}_{s,\varphi}$ with ${\mathcal C}_{i,\varphi}\in EC^{k_i}_{C_i,z}$, $i=1,...,s$.

\begin{lemma}\label{ll5}
In the above notation, the tangent cone $\widehat T_0EC^{\overline k}_{C,z}$
is a linear subspace of $R(C,z)$ of codimension $k+\delta(C,z)=\codim EC^{\overline k}_{C,z}$.
\end{lemma}

{\bf Proof.} It is sufficient to treat the case of an irreducible singularity $(C,z)$.
Let $\varphi$ be a generic element of a component of $EC^k_{C,z}$. The tangent space $T_\varphi EC^k_{C,z}$ at $\varphi$ can be identified with the space
$$\{\psi\in R(C,z)\ :\ \psi(\Sing({\mathcal C}_\varphi))=0,\qquad\qquad\qquad$$
$$\qquad\qquad\qquad\ord\psi\big|_P\ge3\ \text{for each cuspidal local branch}\ P\},$$
and hence the limit of each sequence of tangent spaces $T_\varphi EC^k_{C,z}$ as $\varphi\to0$ is contained in the linear space
$$\{\psi\in R(C,z)\ :\ \ord\psi\big|_{C,z}\ge2\delta(C,z)+k\}\ .$$
It remails to notice that
$$\codim\{\psi\in R(C,z)\ :\ \ord\psi\big|_{C,z}\ge2\delta(C,z)+k\}=\delta(C,z)+k\ .$$
The latter follows, for instance, from \cite[Propositions 5.8.6 and 5.8.7]{CA}.
\proofend

As a corollary we obtain

\begin{proposition}\label{lp6}
Given a real singularity $(C,z)$ splitting into irreducible (over $\C$) irreducible components $(C_i,z)$, $i=1,...,s$, and a sequence $\overline k=(k_1,...,k_s)$ satisfying (\ref{le4}) and an extra condition $k_i=k_j$ as long as $(C_i,z)$ and $C_j,z)$ are complex conjugate, the locus $EC^{\overline k}_{C,z}$ is real, and the number
$W(C,z,EC^{\overline k}_{C,z},L)$ does not depend on the choice of $L$.
\end{proposition}

The proof literally coincides with the proof of Proposition \ref{lp3}.

\section{Example: singularities of type $A_n$}\label{lsec4}

A complex singularity of type $A_n$ is analytically isomorphic to the canonical one $\{y^2-x^{n+1}=0\}\subset(\C^2,0)$, and its miniversal deformation
can be chosen to be
$$\left\{y^2-x^{n+1}-\sum_{i=0}^{n-1}a_ix^i=0\right\}_{a_0,...,a_{n-1}\in(\C,0)}$$ with the base $B(A_n)=\{(a_0,...,a_{n-1})\in(\C^n,0)\}$.

\begin{lemma}\label{ll6} (1) For any $n\ge1$, and $1\le i\le\delta(A_n)=\left[\frac{n+1}{2}\right]$,
\begin{equation}\widehat T_0EG^i_{A_n}=\{a_0=...=a_{i-1}=0\}\subset B(A_n)\label{ne10}\end{equation} the linear subspace of codimension $i=\codim EG^i_{A_n}$.

(2) If $n$ is odd, then $EC_{A_n}=EG_{A_n}$. If $n$ is even, than $$\widehat T_0EC_{A_n}=\{a_0=...=a_k=0\}\subset B(A_n)\ .$$
\end{lemma}

{\bf Proof.} 
Let $(C,z)$ be a canonical singularity of type $A_n$. The tangent space to $EG^i_{C,z}$ at a generic element $\varphi$ consists of $\psi\in B(C,z)$ such that ${\mathcal C}_\psi$ passes through all $i$ nodes of
${\mathcal C}_\varphi$, and hence, $({\mathcal C}_\psi\cdot{\mathcal C}_\varphi)_{D(C,z)}\ge2i$. It follows that the limit of any sequence of these tangent spaces
as $\varphi\to0$ is contained in the linear space $\{\psi\in B(C,z)\ :\ ({\mathcal C}_\psi\cdot C)_z\ge2i\}$, which one can easily identify with the space
in the right-hand side of (\ref{ne10}). So, the first claim of the lemma follows for the dimension reason. The same argument settles the second claim.
\proofend

\begin{proposition}\label{lp5}
For any $n\ge1$ and $k\ge1$, we have
$$\mt EG^i_{A_n}=\binom{n+1-i}{i},\quad\text{for all}\quad i=1,...,\delta(A_n)=\left[\frac{n+1}{2}\right]\ ,$$ $$\text{and}\quad\mt EC(A_{2k})=k\ .$$
\end{proposition}

\begin{remark}\label{nr1}
The multiplicities $\mt EG^i_{A_n}$ were computed in \cite[Section 5, page 540]{She}. Here we provide another, more explicit computation, which will be used below for computing singular Welschinger invariants.
\end{remark}

{\bf Proof.}
(1) If $n+1=2i$, then $EG^i_{A_n}=EG(A_n)=EC(A_n)$ is smooth; hence, the multiplicity equals $1$. Thus, suppose that $n+1>2i$.
By Lemma \ref{ll6}(1), the question on $\mt EG^i_{A_n}$ reduces to the following one: How many polynomials $P(x)$ of degree $\le i-1$ satisfy the condition
\begin{equation}x^{n+1}+x^i+P(x)=Q(x)^2R(x)\ ,\label{ne30}\end{equation} where $Q,R$ are monic polynomials of degree $i$, $n+1-2i$, respectively?

Combining relation (\ref{ne30}) with its derivative, we obtain
$$(n+1-i)x^i+((n+1)P-xP')=\left((n+1)QR-2xQ'R-xQR'\right)Q\ ,$$ which immediately yields
\begin{equation}(n+1)QR-2xQ'R-xQR'=n+1-i\ .\label{ne31}\end{equation}
Substituting
$$Q(x)=x^i+\sum_{j=1}^i\alpha_jx^{i-j},\quad R(x)=x^{n+1-2i}+\sum_{j=1}^{n+1-2i}\beta_jx^{n_1-2i-j}$$ into (\ref{ne31}), we obtain that the terms of the top degree $n+1-i$ cancel out, while the coefficients of $x^m$, $m=0,...,n-i$, yield the system of equations
\begin{equation}\begin{cases}&2\alpha_1+\beta_1=0,\\
&2j\alpha_j+j\beta_j+\sum_{0<m<j}c_{jm}\alpha_{j-m}\beta_m=0,\quad j=2,...,n-i,\\
&(n+1)\alpha_i\beta_{n+1-2i}=n+1-i,\end{cases}\label{ne32}\end{equation}
where we assume $\alpha_j=0$ as $j>i$ and $\beta_j=0$ as $j>n+1-2i$.

Suppose that $\deg Q=i\ge\deg R=n+1-2i$. From the $(n+1-2i)$ first equations in (\ref{ne32}) we express $\beta_j$ as a polynomial in $\alpha_1,...,\alpha_i$ of homogeneity degree $j$, while $\alpha_m$ has weight $m$, for all $j=1,...,n+1-2i$. Substituting these expressions into the other equations, we obtain a system of $i$ equations in $\alpha_1,...,\alpha_i$ of homogeneity degrees $n+2-2i,...,n+1-i$, respectively. Thus, (cf. the computation in
\cite[Section G, Example 1]{FGS})\emph{} the number of solutions (counted with multiplicities) appears to be
$$\frac{(n+2-2i)\cdot...\cdot(n+1-i)}{i!}=\binom{n+1-i}{i}$$ as required. In the same way we treat the case when $\deg Q=i\le\deg R=n+1-2i$.

\smallskip

(2) For $n=2k$, by Lemma \ref{ll6}(2), the question on $\mt EC(A_{2k})$ reduces to the following one: How many polynomials $P(x)$ of degree $k$ satisfy the condition
\begin{equation}x^{2k+1}+x^{k+1}+P(x)=Q(x)^2(x+\beta)^3\ ,\label{ne33}\end{equation} where $Q(x)$ is a monic polynomial of degree $k-1$?

The preceding argument subsequently gives an equation
$$(2k+1)(x+\beta)Q-3xQ-2Q'(x+\beta)=k$$ with $Q(x)=x^{k-1}+\sum_{j=1}^{k-1}\alpha_jx^{k-1-j}$, which develops into the system
\begin{equation}\begin{cases}&2\alpha_1+3\beta=0,\\
&(2j+2)\alpha_j+(2j+3)\alpha_{j-1}\beta=0,\quad j=2,...,k-1,\\
&(2k+1)\alpha_{k-1}\beta=k+1,\end{cases}\label{ne34}\end{equation}
admitting a simplification of the form
$$\alpha_j=\nu_j\beta^j,\ j=1,...,k-1,\quad (2k+1)\nu_{k-1}\beta^k=k+1$$ with some $\nu_1,...,\nu_{k-1}\in\Q$.
So, we finally obtain $k$ solutions
as required.
\proofend

\smallskip

Now we pass to the real setting.
The complex singularity of type $A_n$ has a unique real form $y^2=x^{2k+1}$ if $n=2k$, and has two real forms $y^2=x^{2k}$ and $y^2=-x^{2k}$ (denoted by $A_{2k-1}^h$ and $A_{2k-1}^e$,
respectively) if $n=2k-1$.

\begin{lemma}\label{nl2} (1) For all $k\ge1$ and $i=1,...,k$, there exist singular Welschinger invariants
\begin{equation}W\left(A^h_{2k-1},EG^i_{A^h_{2k-1}}\right),\ W\left(A^e_{2k-1},EG^i_{A^e_{2k-1}}\right),\ \text{and}\ W\left(A_{2k},EG^i_{A_{2k}}\right)\ .\label{ne20}\end{equation}
(2) Furthermore,
$$W^{eg}(A_{2k-1}^e)=(-1)^k,\quad W^{eg}(A_{2k-1}^h)=1\ ,$$
$$W^{eg}(A_{2k})=\begin{cases}0,\quad &k\equiv1\mod2,\\
1,\quad &k\equiv0\mod2,\end{cases}$$ $$W^{ec}(A_{2k})=\begin{cases}0,\quad & k\equiv0\mod2,\\
1,\quad & k\equiv1\mod2.\end{cases}$$
\end{lemma}

{\bf Proof.}
The existence of the invariants (\ref{ne20}) follows from
Lemma \ref{ll6} and the argument used in the proof of Propositions \ref{lp1} and \ref{lp2}.

Since $\mt EG(A_{2k-1})=1$, we have $W^{eg}=\pm1$ for $A^h_{2k-1}$ and $A^e_{2k-1}$. More precisely, an equigeneric nodal deformation of $A^h_{2k-1}$ has the form $y^2-Q(x)^2=0$, $\deg Q=k$, and hence it has only hyperbolic real nodes, i.e., $W^{eg}(A^h_{2k-1})=1$, while an equigeneric nodal deformation of $A^e_{2k-1}$ has the form
$y^2+Q(x)^2=0$, $\deg Q=k$, and hence it has only elliptic real nodes, whose number is of the same parity as $k$, i.e., $W^{eg}(A^e_{2k-1})=(-1)^k$.

Consider singularities $A_{2k}$. For $EG(A_{2k})=EG^k_{A_{2k}}$, system (\ref{ne32}) takes the form
$$\begin{cases}&2\alpha_1+\beta_1=0,\\
&2j\alpha_j+(2j-1)\alpha_{j-1}\beta_1=0,\quad j=2,...,k,\\
&(2k+1)\alpha_k\beta_1=k+1,\end{cases}$$ which yields
$$\alpha_j=\lambda_j\beta_1^j,\ (-1)^j\lambda_j>0,\ j=1,...,k,\quad \lambda_k\beta^{k+1}=\frac{k+1}{2k+1}\ .$$
So, if $k$ is odd, we have no real solutions, and hence $W^{eg}(A_{2k})=0$. If $k$ is even, than we have a unique real solution such that $\beta_1>0$ and
$(-1)^j\alpha_j>0$. That is, $Q(x)$ has only positive real roots (if any), and hence the curve $y^2-(x+\beta_1)Q(x)^2=0$ has only hyperbolic real nodes, i.e.,
$W^{eg}(A_{2k})=1$.

In the same manner we analyze system (\ref{ne34}) and obtain the values of $W^{ec}(A_{2k})$ as stated in the lemma.
\proofend

\begin{remark}\label{nr2}
(1) The problem of computation of the invariants $W^{eg}$ and $W^{ec}$ for arbitrary real singularities (even for quasihomogeneous singularities)
remains widely open. A possible relation to enumerative invariants of (global) plane algebraic curves could be a key to this problem.

(2) The values of $W^{eg}$ and $W^{ec}$ for $A_n$-singularities are $0$ or $\pm1$. The same can be showed for other simple singularities. Is it true
for an arbitrary real singularity?
\end{remark}

\end{document}